\documentclass{article}
\usepackage[utf8]{inputenc}
\usepackage{titling}
\usepackage{lipsum}
\usepackage[space]{grffile}
\usepackage{graphicx}
\usepackage{url}
\usepackage{amsmath}
\usepackage{amssymb}
\usepackage{upgreek}

\newcommand{\norm}[1]{\left\lVert#1\right\rVert}

\usepackage[thmmarks, amsmath, thref]{ntheorem}

\newtheorem{observation}{Observation}[section]
\newtheorem{theorem}{Theorem}[section]

\newtheorem{lemma}[theorem]{Lemma}
\newtheorem{proof}{Proof}[theorem]

\newtheorem{definition}{Definition}[section]
\newtheorem{example}{Example}[section]

\title{\bf Canonical Tensor Scaling}
\author{{\bf Tung D.\ Nguyen} and {\bf Jeffrey Uhlmann}\vspace{4pt} \\
University of Missouri - Columbia}

\begin{document}
\maketitle

\section{Introduction}

A 1992 algorithm by Rothblum \& Zenios \cite{RZalgorithm}, which we will refer to as the RZ algorithm, 
computes a positive diagonal scaling of an arbitrary matrix so that the products of 
the nonzero normed elements of each row and column (i.e., zero elements are ignored) 
are equal to user-specified constants. It has recently been recognized that
in the case in which the constants are all unity, a unique scaling always exists and
can be used as a canonical scaling for performing scale-consistent or scale-invariant operations on matrices \cite{UC18}. For example, the algorithm can be applied to define a generalized matrix inverse that is consistent with respect to diagonal scalings (e.g., due to changes of units on variables), as opposed to the Moore-Penrose pseudoinverse which provides consistency with respect to unitary transformations (e.g., rotations) \cite{UC18}. Scale consistency is critical to applications in which performance is expected to be invariant with respect to the choice of state variables, e.g., whether lengths are represented in centimeters, meters, or kilometers in robotics applications \cite{rga,BoZ1,BoZ2,JKU18} or unknown units specific to each individual's ratings of products in the context of recommender systems. 

In this paper we generalize the canonical positive scaling of rows and columns of a matrix to the
scaling of selected-rank subtensors of an arbitrary tensor. We expect our results and framework will prove useful for sparse-tensor completion required for generalizations of the recommender system problem beyond a matrix of user-product ratings to multidimensional arrays involving coordinates based both on user attributes (e.g., age, gender, geographical location, etc.) and product/item attributes (e.g., price, size, weight, etc.).

\section{Definitions}

The objective of this paper is to identify a generalized canonical scaling of subtensors of a given tensor. In the case of a matrix, i.e., a 2-tensor, the rows and columns can be interpreted as the 1-dimensional subtensors that can be scaled by the RZ algorithm. In order to generalize beyond matrices to arbitrary $d$-dimensional tensors, or $d$-tensors, we must define notation for operating on the $k$-dimensional subtensors that are to be scaled, where $1\leq k < d$ for integer parameter $k$. This requires distinguishing the overlapping $k$-dimensional subtensors with respect to each of the $d$ coordinates. In the $d=2$ case of a matrix with $k = 1$, for example, the dimensional label set $\{1,2\}$ can be interpreted with $(1)$ denoting the set of row vectors and $(2)$ representing the set of column vectors. In other words, the tuple $(1,i)$ designates the $i$th row vector while $(2,j)$ designates the $j$th column vector\footnote{Matrices are special in that the terms {\em rows} and {\em columns} have been conventionally defined to distinguish two sets of overlapping $1$-dimensional subtensors/vectors. For our purposes the particular choice of the ordering of the labels, i.e., the distinction between a given matrix and its transpose, is irrelevant. More generally, there is no natural ordering of the $d$ orientations of a $d$-tensor with names analogous to {\em rows} and {\em columns}, so for notational purposes any arbitrary choice will work.}, and $(2,i,j)$ would specify the $j$th element of the $i$th column.  

In the $d=3$ case of a  $3$-tensor, the choice $k\in \{1,2\}$ can be chosen to specify whether our focus is on $2$-dimensional subtensors (matrices) or $1$-dimensional subtensors (vectors). For example, $k=1$ defines $d=3$ sets of $1$-dimensional subtensors, or fibers, each set of which can be thought of intuitively as containing all of the vectors parallel to one of the $3$ coordinate axes, e.g., set $(1)$ refers to all $1$-dimensional subtensors that live in $\mathbb{R}^{n_1\times 1\times 1} \equiv \mathbb{R}^{n_1}$, and $k=2$ defines all possible $2$-dimensional subtensors, or {\em slices}, comprising $\binom{3}{2}= 3$ sets of slices, the span of which is determined by one of the pairs of indices $(1,2)$, $(2,3)$, or $(1,3)$. For instance, a pair $(1,3)$ spans all possible $2$-dimensional subtensors that live in $\mathbb{R}^{n_1\times 1 \times n_3} \equiv \mathbb{R}^{n_1 \times n_3}$. In other words, it is possible to relabel these pairs as sets $(1)$, $(2)$, and $\left( \binom{3}{2}\right) = 3.$ An important property of the subtensors in each set is that they are disjoint from each other and can be stacked to form the original $d$-dimensional tensor. In general, the $k$-dimensional subtensors of the $d$-dimensional tensor can be partitioned into $\left(\binom{d}{k}\right)$ sets, corresponding to the number of $k$-tuples chosen from $d$ dimensional indices. An alternative but equivalent interpretation is that there are $\binom{d}{k}$ ways of dividing a $d$-dimensional tensor into $k$-subtensors such that all subtensors in each set are disjoint and any two subtensors in two different sets are different. This leads to the following definition. 

\begin{definition} - $\textup{\bf Sets of subtensors:}$ Let $V_i = \mathbb{R}^{n_i}$ and assume a $d$-dimensional tensor $A \in V_1 \times V_2 \times \cdots \times V_d$ and a positive integer $1 \leq k < d$. From the set of dimensional indices $[d] = \{1,\dots, d\}$ we can obtain all $\binom{d}{k}$ possible tuples of $k$ dimensional indices. Specifically, for any tuple $\{n_{i_1}, \cdots, n_{i_k}\},$ we enumerate all possible $k$-dimensional subtensors of $A$ in dimensions $V_{i_1} \times \cdots \times V_{i_{k}}$ to obtain the set of subtensors $(i)$. Altogether we obtain $\binom{d}{k}$ sets of $k$-dimensional subtensors relabeled as $(1), \dots, \left( \binom{d}{k}\right)$. We denote the number of subtensors in a set $(i)$ as its cardinality $|(i)|$ and label its elements as
\begin{equation}
    (i) = \left\{A_{s,i} \mid 1 \leq s \leq |(i)| \right\}
\end{equation}
for $1 \leq i \leq \binom{d}{k}$.
\end{definition}
By this definition, all subtensors in each set $(i)$ are disjoint and can combine to reconstruct the tensor $A$. Each set $(i)$ corresponds to a partition of the $d$-dimensional tensor into mutually disjoint $k$-dimensional subtensors. This property implies the following observation.

\begin{observation} $\textup{\bf Tensor position in subtensors:}$ Given a $d$-dimensional tensor $A$, an element at position $\Vec{\alpha} = (\alpha_1, \alpha_2, \cdots, \alpha_d) $ of a tensor $A \in V_1 \times V_2 \times \cdots \times V_d$  is defined as $A(\Vec{\alpha}) \equiv A(\alpha_1, \cdots, \alpha_d).$ Then there is a unique subtensor with index $(s_i,i)$ from a set $(i)$, for $1 \leq s_i \leq |(i)|$ and $1 \leq i \leq \binom{d}{k}$ such that element $A(\Vec{\alpha})$ is contained in the subtensor $A_{s_i,i}$, denoted as $A(\Vec{\alpha}) \sqsubset A_{s_i,i}$. In other words, there is a unique list of subtensors $\left\{A_{s_1,1}, \cdots, A_{s_i, i}, \cdots, A_{s_{\binom{d}{k}}, \binom{d}{k}} \right\}$ for each $A(\Vec{\alpha})$ such that for each of its element $A_{s_j,j},$ $A(\Vec{\alpha}) \sqsubset A_{s_j,j}$.
\label{Position_rep}
\end{observation}
Throughout, we consider two conformant $d$-dimensional tensors $A$ and $A'$. We establish a similarity (correspondence) between each partition of $\binom{d}{k}$ sets of $k$-dimensional subtensors of $A$ and $A'$ as follows.

\begin{definition} -  $\textup{\bf Similar sets of subtensors:}$
Given two tensors $A$ and $A'$ that are both in $V_1 \times V_2 \times \cdots \times V_d$. Each of them has $\binom{d}{k}$ sets of $k$-dimensional subtensors labeled as $(1), \dots, \left(\binom{d}{k}\right)$, but each k-dimensional subtensor of $A'$ corresponding to a pair $(s_i, i)$ labels as $A'_{s_i,i}$ and that of $A$ is $A_{s_i,i}.$ Then we say both $\binom{d}{k}$ sets of $k$-dimensional subtensors of $A$ and $A'$ are \textbf{similar} if for each pair $A'_{s_i, i}$ and $A_{s_i,i}$, the following condition satisfies for every possible position $\Vec{\alpha} = (\alpha_1, \alpha_2, \cdots, \alpha_d) $,
\begin{equation}
    A'(\Vec{\alpha}) \sqsubset A'_{s_i,i}  \Leftrightarrow A(\Vec{\alpha}) \sqsubset  A_{s_i,i}
\end{equation}
\label{similar_subtensors} \end{definition} 
Now we are equipped to define the scaling process. For notational convenience we define a list $M_{(i)}$ of elements such that $c M_{(i)}$ multiplies/scales each element of every $k$-dim subtensor of $|(i)|$ as a single scaling process, and the product of a scalar in $M_{(i)}$ and a subtensor $A_{s_i,i}$ scales every elements of that subtensor and thus scales the tensor itself. We now formally define both the set $M_{(i)}$ and notation for this scaling operator.

\begin{definition} - $\textup{\bf Scaling a tensor by a set of subtensors:}$ 
Given a tensor $A \in V_1 \times V_2 \times \cdots \times V_d$, the list of scaling indices is defined as
$$M_{(i)} = \{M_{s,i} \in \mathbb{R}_{\neq 0} \mid 1 \leq s \leq |(i)|\},$$
and the scaling process for the set $(i)$ of tensor $A$ by tensor $A'$, which has the same dimension as tensor $A$. Consider partitions into $\binom{d}{k}$ sets of $k$-dimensional subtensors of both tensors $A$ and $A'$ such that the partitions are \textbf{similar} by definition \ref{similar_subtensors}. The scaling operation $\circ_{(i)}$ is defined as
\begin{equation}
    A' \equiv A\circ_{(i)}M_{(i)}
\end{equation}
such that the Hadamard (elementwise) product of elements of $A$ corresponding to $A_{s,i}$ and $M_{s,i}$ results in the scaling of elements of $A'$ corresponding to $A'_{s,i}$ as
\begin{equation}
    A'_{s, i} = A_{s,i} \circ M_{s,i} 
\end{equation}
for $1 \leq s \leq |(i)|$. 
\label{scalingbyHadamard}
\end{definition}
\begin{example}
Consider the case when $d = 3$ and $k=2$. Let $A \in \mathbb{R}^{3\times 4\times 2}$ be a tensor with the following 2 frontal slices.
\begin{equation}
 A_1 = \begin{bmatrix}
    1 & 4 & 7 & 10 \\
    2 & 5 & 8 & 11 \\
    3 & 6 & 9 & 12
\end{bmatrix}
, A_2 =  
\begin{bmatrix}
    13 & 16 & 19 & 22 \\
    14 & 17 & 20 & 23 \\
    15 & 18 & 21 & 24
\end{bmatrix} 
\end{equation}
Then the set $(3)$ contains the subtensors $A_1$ and $A_2$ and the scaling list $M_{(3)} = \{\alpha_1, \alpha_2\},$ then the scaling process by a list $M_{(3)}$ gives the new tensor $S$ as
\begin{equation}
 S_1 = \begin{bmatrix}
    \alpha_{1} & 4\alpha_{1} & 7\alpha_{1} & 10\alpha_{1} \\
    2\alpha_{1} & 5\alpha_{1} & 8\alpha_{1} & 11\alpha_{1} \\
    3\alpha_{1} & 6\alpha_{1} & 9\alpha_{1} & 12\alpha_{1}
\end{bmatrix}
, S_2 =  
\begin{bmatrix}
    13\alpha_{2} & 16\alpha_{2} & 19\alpha_{2} & 22\alpha_{2} \\
    14\alpha_{2} & 17\alpha_{2} & 20\alpha_{2} & 23\alpha_{2} \\
    15\alpha_{2} & 18\alpha_{2} & 21\alpha_{2} & 24\alpha_{2}
\end{bmatrix} 
\end{equation}

\end{example}
At this point we have defined a particular type of structured scaling that can be applied to elements of specified sets of subtensors of a given tensor. In the next section we iteratively apply scalings of this kind to obtain a unique scaling with prescribed properties.

\section{Unique Tensor Scaling}
Using the representation and notation defined in the previous section, we demonstrate the uniqueness and correctness of our generalization of the RZ algorithm for matrices (i.e., $n=2$ and $k=1$) to arbitrary tensors with $1\leq k<n$.  

\subsection{Defining the $k$-dim scaling program:}
\subsubsection{Program I:}
We now follow and extend the approach of Rothblum and Zenios (RZ) in \cite{RZalgorithm} for the iterative scaling of 1-dimensional row and column vector sets of a matrix to the general case for arbitrary scaling of $\binom{d}{k}$ sets of subtensors of a given nonnegative tensor $A$, which may be obtained by replacing the elements of a given tensor $A'$ with their magnitudes. The algorithm begins with tensor $A \in V_1 \times V_2 \times \cdots \times V_d$ and a list of strictly positive numbers $S = \{S_{s_i,i} \mid 1 \leq s_i \leq |(i)|$\text{ and }$1 \leq i \leq \binom{d}{k} \}$. Additionally, we can simplify the representation of the tensor scaling problem by the following definitions:
\begin{definition} - \textup{\bf Subtensors by an element}
Given a unique list of subtensors $\{A_{s_1,1}, \cdots, A_{s_i, i},
\cdots,\\ A_{s_{\binom{d}{k}}, \binom{d}{k}} \}$ such that each of them contains $A(\Vec{\alpha})$, we denote the set of all pairs $(s_i, i)$ for each $A(\Vec{\alpha})$ as $\upsigma_A(\Vec{\alpha}) = \left\{(s_1, 1), \dots, \left(s_{\binom{d}{k}}, \binom{d}{k}\right)\right\} = \left\{(s_i, i) \mid 1\leq i \leq \binom{d}{k}, A(\Vec{\alpha}) \sqsubset A_{s_i,i} \right\}$.
\label{s_i_get_tensor_ele}
\end{definition}
\begin{definition} - \textup{\bf Support set of a subtensor} Define the support set of subtensor $A_{s_i,i}$ as 
\begin{equation}
    \upsigma_{A}(s_i,i) = \left\{ \Vec{\alpha} \mid (s_i,i) \in \upsigma_A(\Vec{\alpha}), A(\Vec{\alpha}) \neq 0 \right\}
\end{equation}
where their union is the set of nonzero elements in the tensor $A$, denoted by $\upsigma(A)$,
\begin{equation}
    \upsigma(A) = \bigcup_{(s_i,i)} \upsigma_A(s_i,i).
\end{equation} \label{support_set}
\end{definition}
Let the scaling elements of $d\choose k$ list $M_{(i)}$ of scaling elements be applied in the scaling process of tensor $A$ in label order from set $(1)$ to set $\binom{d}{k}$. We observe that the choice of ordering is free due to the the commutative property of the scaling process.
\begin{observation} The scaling of indices $i_1, i_2 \in \left[\binom{d}{k}\right]$ is commutative:
\begin{equation}
    A \circ_{(i_1)}  M_{(i_1)} \circ_{(i_2)}  M_{(i_2)} = A \circ_{(i_2)}  M_{(i_2)} \circ_{(i_1)}  M_{(i_1)}.
\end{equation}
\end{observation}
In other words, a permutation of the dimensional indices will not affect the correctness of the following problem. Using definitions \ref{scalingbyHadamard}, \ref{s_i_get_tensor_ele}, and \ref{support_set}, we define the scaling problem:
\\\\
\textbf{Problem Statement}: Find a positive tensor tensor $A' \in V_1 \times V_2 \times \cdots \times V_d$ and $\binom{d}{k}$ scaling lists  $M_{(i)}$ of positive elements with cardinality $|(i)|$ such that:
\begin{equation}
    A' = A \circ_{(1)}  M_{(1)} \cdots \circ_{\left( \binom{d}{k}\right)}  M_{\left( \binom{d}{k}\right)} 
\end{equation}
or 
\begin{eqnarray*}
a'(\Vec{\alpha}) \equiv   \left\{ 
\begin{array}{c}
a(\Vec{\alpha})\cdot \prod\limits_{(s_i, i) \in \upsigma_A(\Vec{\alpha})} m_{s_i, i}  \qquad \text{if} \qquad  \vec{\alpha} \in \upsigma(A) \\ 
0 \quad \text{  otherwise  } 
\end{array}
\right.
\end{eqnarray*}
and for each subtensor by a pair $(s_i, i)$
\begin{equation}
\prod_{\Vec{\alpha} \in \upsigma_{A}(s_i,i)} A'(\Vec{\alpha})  = S_{s_i, i}
\end{equation}
A logarithm conversion for tensor $A, A',$ and $S$ is applied via the below-defined operator $\Hat{L}$ to obtain tensors $a, a'$, and $s$:
\begin{equation}
   \begin{aligned}
   &a(\Vec{\alpha}) \equiv \Hat{L}(A(\vec{\alpha})) \equiv   \left\{ 
\begin{aligned}
&ln(A(\Vec{\alpha})) \qquad \text{if} \qquad \Vec{\alpha} \in \upsigma(A)  \\ 
&0 \quad \text{  otherwise  }
\end{aligned}
\right.
\\
    &a'(\Vec{\alpha}) \equiv \Hat{L}(A'(\vec{\alpha}))
   \\
   &s_{s_i, i} \equiv ln(S_{s_i, i})
   \\
   &m_{s_i, i} \equiv ln(M_{s_i, i})
   \end{aligned}
\end{equation}
Having performed the logarithm conversion of the operator $\Hat{L}$, we present the following program:
\\
\\
\textbf{Program I}: Finding the tensor $a' \in V_1 \times V_2 \times \cdots \times V_d$ and $d\choose k$ lists $M_{(i)}$ such that:
\begin{eqnarray*}
a'(\Vec{\alpha}) \equiv   \left\{ 
\begin{array}{c}
a(\Vec{\alpha}) + \sum\limits_{(s_i, i) \in \upsigma_A(\Vec{\alpha})} m_{s_i, i}  \qquad \text{if} \qquad  \vec{\alpha} \in \upsigma(A) \\ 
0 \quad \text{  otherwise  } 
\end{array}
\right.
\end{eqnarray*}
and for each subtensor by a pair $(s_i, i)$
\begin{equation}
\sum_{\Vec{\alpha} \in \upsigma_{A}(s_i,i)} a'(\Vec{\alpha})  = s_{s_i, i}
\end{equation}

\subsubsection{Program II:}
The RZ model provides a basis for establishing the existence and uniqueness properties of our tensor generalization. Briefly, let $a \in R^p, b \in R^q, $ and $C\in R^{p\times q}$ be the original convex optimization problem of finding a vector $a' \in R^p$ and $\omega \in R^q$ such that
\begin{equation}
    a'^{T}=a^T + \omega^T C
\end{equation}
and 
\begin{equation}
    C a' = b
\end{equation}
It is proven in \cite{RZalgorithm} that this problem is equivalent to the following optimization problem, for which the properties of uniqueness and existence can be established: 
\begin{equation}
\begin{aligned}
    \text{min } 2^{-1} &\sum^p_{j=1} (x_j-a_j)^2 & 
    \\
    & \text{subject to} & Cx = b.
\end{aligned}
\end{equation}
To establish the corresponding properties of the RZ algorithm for our \textbf{Program I}, we apply the following transformations. We begin with a transformation to {\em unfold} tensors $a$ and $a'$ according to the following definition.
\begin{definition} - $\textup{\bf Unfolding tensor:}$ Given a tensor $A = V_1 \times V_2 \times \cdots \times V_d$ and $P = \prod^d_{j=1} n_j$. The unfolding vector $A_P \in \mathbb{R}^{P}$ defined via the following mapping: 
$$A (\Vec{\alpha}) \mapsto A_P(J(\Vec{\alpha}))$$ for 
\begin{equation}
\begin{aligned}
    \nonumber J(\Vec{\alpha}) &= 1 + \sum\limits_{s=1}^d \left((\alpha_s - 1)\prod\limits_{m=1}^{s-1}
    n_m\right)
    \\
    &= \alpha_1 + (\alpha_2 -1)n_1 +  \cdots + (\alpha_d-1)n_1\cdots n_{d-1}
\end{aligned}{}
\end{equation}
\label{gettheorderforC}
\end{definition}
Regarding $Cx =b$ and $\omega^T C$, we show (following \cite{RZalgorithm}) by the Kuhn-Tucker (KT) conditions that  \textbf{Program I} is equivalent to
\\\\
\textbf{Program II}: For any pair $(s_i, i)$: 
\begin{equation}
\begin{aligned}
 \text{find min } & 2^{-1}\sum\limits_{\Vec{\alpha} \in \upsigma_{A}(s_i,i)}(x\left(\Vec{\alpha}\right)  - a\left(\Vec{\alpha}\right))^2, 
 \\
 & \text{subject to}  \sum\limits_{\Vec{\alpha} \in \upsigma_{A}(s_i,i)} x(\Vec{\alpha})  = s_{s_i, i}
\end{aligned}
\end{equation}
where in a slight abuse of notation we have used the tensor forms and unfolding vector forms of $a$, $a'$, and $x$ interchangeably.

\subsubsection{Proof:}

We define the node-arc incidence matrix of tensor $A$ as a matrix $C$ of real numbers with the objective of transforming the 
tensor problem to a matrix problem for which the solution becomes equivalent to the RZ solution for $C$. Specifically, we establish formulas such that each position of $\omega^T C$ defined by a function $J\left(\Vec{\alpha}\right)$ equals the sum of $m_{s_i,i}$ of \textbf{Program I} and each position on $Cx$ defined by a pair $(s_i,i)$ equals to the nonzero sum of $x$ of \textbf{Program II}.  The basic idea is to use the unfolding tensor formula from definition \ref{gettheorderforC}. 
We provide the detailed construction in \textbf{Appendix 1}, and in brief the vector $\omega$ and $C$ matrix satisfy 
\begin{equation}
    [\omega^T C]_{J(\Vec{\alpha})} = \sum\limits_{(s_i, i) \in \upsigma_A(\Vec{\alpha})} m_{s_i, i}.  
\end{equation}
where the matrix $b$ is defined using $s_{s_i, i}$ such that $Cx = b$ is equivalent to

\begin{equation}
  \sum\limits_{\Vec{\alpha} \in \upsigma_{A}(s_i,i)} x(\Vec{\alpha})  = s_{s_i, i}.  
\end{equation}
Satisfaction of the KT conditions guarantees
\begin{equation}
    \begin{aligned}
    0 = & \dfrac{\partial}{\partial x(\Vec{\alpha})} \left[ 2^{-1} \sum\limits_{\Vec{\alpha} \in \upsigma_{A}(s_i,i)}(x\left(\Vec{\alpha}\right)  - a\left(\Vec{\alpha}\right))^2 - \omega^T (Cx - b)\right] \\
    = & x(\Vec{\alpha}) - a(\Vec{\alpha}) - [\omega^T C]_{J(\Vec{\alpha})}  
    = x(\Vec{\alpha}) - a(\Vec{\alpha}) - \sum\limits_{(s_i, i) \in \upsigma_A(\Vec{\alpha})} m_{s_i, i}.
    \end{aligned}
\end{equation}
We now exploit the equivalence of our transformed formulaton to that of the optimization problem of section 3 of \cite{RZalgorithm} to establish the following:
\begin{theorem}{(Characterization)} The following statements are equivalent:
\begin{enumerate}
    \item There exists a solution to \textbf{Program I}
    \item There exists a solution to \textbf{Program II}
    \item There exists a tensor $A' \in V_1 \times V_2 \times \cdots \times V_d$ that satisfies \textbf{Program I} and has $\upsigma(A) = \upsigma(A')$
    \item \textbf{Program II} is feasible 
    \item \textbf{Program II} has a optimal solution
    \item Given $\binom{d}{k}$ list of real numbers $\mu_{(i)} = \{\mu_{s_i,i} \in \mathbb{R} \mid 1 \leq s_i \leq |(i)|\}$ with cardinality $|(i)|$ satisfying
    $\sum\limits_{(s_i,i) \in \sigma_A(\Vec{\alpha})} \mu_{s_i,i} = 0$ for all 
    $\Vec{\alpha} \in \sigma(A)$, then
    
    \begin{equation}
     \prod\limits_{i=1}^{\binom{d}{k}} \prod\limits_{s_i = 1}^{|(i)|} \left(S_{s_i,i}\right)^{\mu_{s_i,i}} = 1.
    \end{equation}

    \item Similarly from condition 6, 
    \begin{equation} \sum\limits_{i=1}^{\binom{d}{k}} \sum\limits_{s_i = 1}^{|(i)|} s_{s_i,i}\cdot\mu_{s_i,i} = 0 
    \end{equation}

\end{enumerate}{}
 \label{characterize}
\end{theorem}{}
The detailed proof is given in \textbf{Appendix 2}, with Theorem \ref{characterize} playing a fundamental role. Specifically, we know that if there exists a solution $x$ such that $Ax = b$, and vector $\omega$ such that $A\omega = 0$, then $x + \omega$ is another solution. Thus, additional solutions can be obtianed from the known solution of \textbf{Program I} based on statements $(6)$ and $(7)$ in Theorem \ref{characterize}.

\begin{theorem}{(Uniqueness)} There exists at most one tensor $A' \in V_1 \times V_2 \times \cdots \times V_d$ such that there exist $\binom{d}{k}$ scaling lists  $M_{(i)} = \left\{M_{s_i, i} \mid 1 \leq s_i \leq |(i)| \right\}$ of positive elements with cardinality $|(i)|$ such that the solution $\left(A',\left\{ M_{(i)} \mid 1\leq i \leq \binom{d}{k} \right\}\right)$ satisfies \textbf{Program I}. Furthermore, if $\left(A',\left\{ M_{(i)}^1 \mid 1\leq i \leq \binom{d}{k} \right\} \right)$ satisfies \textbf{Program I}, then for $\binom{d}{k}$ lists $T_{(i)}=\{T_{s_i,i} \in \mathbb{R}_{>0}| \mid 1 \leq s_i \leq |(i)|\}$ such that
\begin{equation}
    \prod\limits_{(s_i, i) \in \upsigma_A(\Vec{\alpha})}  T_{s_i, i}  = 1 \qquad \forall \Vec{\alpha}  \in \upsigma(A),
\end{equation} its general solution is 
\begin{equation}
    \left(A',\left\{ MT_{(i)}^1  \mid 1\leq i \leq \binom{d}{k} \right\}\right)
\end{equation}
where $MT_{(i)}^1 = \left\{ M_{s_i, i}^1 T_{s_i, i} \mid 1 \leq s_i \leq |(i)|, T_{s_i, i} \in T_{(i)} \right\}$.
\end{theorem}{}
\begin{proof}
The proof of uniqueness follows from the RZ model, and the general solution follows directly from Theorem \ref{characterize} for the transformed expression of \textbf{Program I}.
\end{proof}

\subsection{Demonstration of $(d-1)-$dim scaling:}
To illustrate, we consider the case of $k=d-1$:
\\\\
\textbf{Problem}: Find a positive tensor tensor $A' \in V_1 \times V_2 \times \cdots \times V_d$ and $d$ scaling lists  $M_{(i)}$ of positive elements with cardinality $|(i)|$ such that:
\begin{equation}
    A' = A \circ_{(1)}  M_{(1)} \cdots \circ_{\left( d\right)}  M_{\left( d\right)} 
\end{equation}
or 
\begin{eqnarray*}
a'(\Vec{\alpha}) \equiv   \left\{ 
\begin{array}{c}
a(\Vec{\alpha})\cdot \prod\limits_{(s_i, i) \in \upsigma_A(\Vec{\alpha})} m_{s_i, i}  \qquad \text{if} \qquad  \vec{\alpha} \in \upsigma(A) \\ 
0 \quad \text{  otherwise  } 
\end{array}
\right.
\end{eqnarray*}
and for each subtensor by a pair $(s_i, i)$ for $1 \leq i\leq d$ 
\begin{equation}
\prod_{\Vec{\alpha} \in \upsigma_{A}(s_i,i)} A'(\Vec{\alpha})  = S_{s_i, i}
\end{equation}
Using operator $\hat{L}$, the logarithmic form of \textbf{Problem} is,
\\\\
\textbf{Program I}: Finding the tensor $a' \in V_1 \times V_2 \times \cdots \times V_d$ and $d$ lists $M_{(i)}$ such that:
\begin{eqnarray*}
a'(\Vec{\alpha}) \equiv   \left\{ 
\begin{array}{c}
a(\Vec{\alpha}) + \sum\limits_{(s_i, i) \in \upsigma_A(\Vec{\alpha})} m_{s_i, i}  \qquad \text{if} \qquad  \vec{\alpha} \in \upsigma(A) \\ 
0 \quad \text{  otherwise  } 
\end{array}
\right.
\end{eqnarray*}
and for each subtensor by a pair $(s_i, i)$ for $1 \leq i\leq d$ 
\begin{equation}
\sum_{\Vec{\alpha} \in \upsigma_{A}(s_i,i)} a'(\Vec{\alpha})  = s_{s_i, i}
\end{equation}
\\
From our framework, this program is equivalent to the following program,
\\\\
\textbf{Program II}: For any pair $(s_i,i)$,
\begin{equation}
\begin{aligned}
 \text{find min } & 2^{-1}\sum\limits_{\Vec{\alpha} \in \upsigma_{A}(s_i,i)}(x\left(\Vec{\alpha}\right)  - a\left(\Vec{\alpha}\right))^2, 
 \\
 & \text{subject to}  \sum\limits_{\Vec{\alpha} \in \upsigma_{A}(s_i,i)} x(\Vec{\alpha})  = s_{s_i, i}
\end{aligned}
\end{equation}
Hence the theorem for the general solution when $k=d-1$ states,

\begin{theorem}{(Uniqueness)}
There exists at most one tensor $A' \in V_1 \times V_2 \times \cdots \times V_d$ such that there exist and $d$ scaling lists  $M_{(i)} = \left\{M_{s_i, i} \mid 1 \leq s_i \leq |(i)| \right\}$ of positive elements with cardinality $|(i)|$ such that the solution $\left(A',\left\{ M_{(i)} \mid 1\leq i \leq d \right\}\right)$ satisfies \textbf{Program I}. Furthermore, if $\left(A',\left\{ M_{(i)}^1 \mid 1\leq i \leq d \right\} \right)$ satisfies \textbf{Program I}, then for $d$ lists $T_{(i)}=\{T_{s_i,i} \in \mathbb{R}_{>0}| \mid 1 \leq s_i \leq |(i)|\}$ such that
\begin{equation}
    \prod\limits_{(s_i, i) \in \upsigma_A(\Vec{\alpha})}  T_{s_i, i}  = 1 \qquad \forall \Vec{\alpha}  \in \upsigma(A),
\end{equation} its general solution is 
\begin{equation}
    \left(A',\left\{ MT_{(i)}^1  \mid 1\leq i \leq d \right\}\right)
\end{equation}
where $MT_{(i)}^1 = \left\{ M_{s_i, i}^1 T_{s_i, i} \mid 1 \leq s_i \leq |(i)|, T_{s_i, i} \in T_{(i)} \right\}$.
\end{theorem}

\section{Algorithm}
To simplify the algorithmic representation, we define a function that maps the number of elements in $A_{s,i}$ from an input as a pair $(s, i)$. This establishes a relationship between the number of nonzero elements in $A_{s_i,i}$ and the cardinality of $\upsigma_A(s_i,i)$
\begin{definition} - $\textup{\bf Number of elements in a subtensor}$
For a set $(i)$, let $\phi$ be a function such that $\phi(s,i)$ is the number of nonzero elements in $A_{s,i}$, for $1 \leq s \leq |(i)|$ and $1 \leq i \leq \binom{d}{k}$.
\end{definition}
One important property is that the cardinality of $\upsigma_A(s_i,i)$ then equals $\phi(s_i,i)$. We establish the following algorithm.
\subsection{The case $k = d-1$:}
Base on the framework of RZ's model in \cite{RZalgorithm}, we establish the following algorithm for case $k=d-1$.
\\ 
\begin{center}
    SCALING ALGORITHM FOR PROGRAM II
\end{center}
\begin{enumerate}
    \item{Step 0 [initialization]:} Set $count \leftarrow 0.$ Select $d$ list of scaling $m_{(i)} = \{m_{s,i} \in \mathbb{R} \mid 1 \leq s \leq |(i)| \}$ and $d$ sets of subtensors $(i)$  by definition \ref{Position_rep} where $d$ is the dimensionality of tensor $A$. Define tensor $a^0 \in V_1 \times V_2 \times \cdots \times V_d$ such that \begin{equation}
     a^0(\Vec{\alpha}) \equiv  \left\{ \begin{aligned}
   & a(\Vec{\alpha}) + \sum\limits_{(s_i, i) \in \upsigma_A(\Vec{\alpha})} m_{s_i, i} \quad & \text{  for} \quad \Vec{\alpha} \in \upsigma(A)  \\  
&0 \quad & \text{  for} \quad \Vec{\alpha} \notin \upsigma(A)
\end{aligned}\right.
\end{equation}

\item{Step $1$ [iterative step over constraint]:} Let $i = 1$ and our goal is to iterate until $i =d.$ For $s_i \in \left[ |(i)|\right],$ let
\begin{equation}
\begin{aligned}
    \rho_{s_i} &= \left[ \phi(s_i,i)\right]^{-1} \left[ s_{s_i,i} -
    \sum\limits_{\Vec{\alpha} \in \upsigma_{A}(s_i,i)} a^{step}(\Vec{\alpha}) 
    \right] 
    \\
    a^{step}(\Vec{\alpha}) & \leftarrow a^{step}(\Vec{\alpha}) + \rho_{s_i} \quad \text{for} \quad \Vec{\alpha} \in \upsigma_{A}(s_i,i)
    \\
    m_{s_i,i}&\leftarrow m_{s_i,i} + \rho_{s_i}
\end{aligned}
\end{equation}
    We then set $i \leftarrow i+1$ and continue until $i = d.$ 
    \item{Step 2:} Set $count \leftarrow count+1$ and return to step 1. We also set the convergence condition as the variance of the tensor elements projected by each of its dimensions is smaller than a chosen small positive $\epsilon$ (though alternative termination criteria may be applied).
    
\end{enumerate}
The time complexity is $ O\left(\sum\limits_{i=1}^{d} \sum\limits_{s_i = 1}^{|(i)|}\phi(s_i,i) \right) \equiv O\left(d\norm{\upsigma(A)}\right)$. This complexity is unaffected by the final step of converting our solution over logarithms to our desired solution as $A(\vec{\alpha}) = \exp(a(\Vec{\alpha}))$, which is proportional to the number of nonzero elements of $A$.

\subsection{Algorithm for general $k < d$:}

\begin{center}
    SCALING ALGORITHM FOR PROGRAM II
\end{center}
\begin{enumerate}
    \item{Step 0 [initialization]:} Set $count \leftarrow 0.$ Select $\binom{d}{k}$ list of scaling $m_{(i)} = \{m_{s,i} \in \mathbb{R} \mid 1 \leq s \leq |(i)| \}$ and $\binom{d}{k}$ sets of subtensors $(i)$  by definition \ref{Position_rep} for $d$ as the dimension of tensor $A$. Define tensor $a^0 \in V_1 \times V_2 \times \cdots \times V_d$ such that \begin{equation}
     a^0(\Vec{\alpha}) \equiv  \left\{ \begin{aligned}
   & a(\Vec{\alpha}) + \sum\limits_{(s_i, i) \in \upsigma_A(\Vec{\alpha})} m_{s_i, i} \quad & \text{  for} \quad \Vec{\alpha} \in \upsigma(A)  \\  
&0 \quad & \text{  for} \quad \Vec{\alpha} \notin \upsigma(A)
\end{aligned}\right.
\end{equation}

\item{Step $1$ [iterative step over constraint]:} Let $i = 1$ and our goal is to iterate until $i =\binom{d}{k}.$ For $s_i \in \left[ |(i)|\right],$ let
\begin{equation}
\begin{aligned}
    \rho_{s_i} &= \left[ \phi(s_i,i)\right]^{-1} \left[ s_{s_i,i} -
    \sum\limits_{\Vec{\alpha} \in \upsigma_{A}(s_i,i)} a^{step}(\Vec{\alpha}) 
    \right] 
    \\
    a^{step}(\Vec{\alpha}) & \leftarrow a^{step}(\Vec{\alpha}) + \rho_{s_i} \quad \text{for} \quad \Vec{\alpha} \in \upsigma_{A}(s_i,i)
    \\
    m_{s_i,i}&\leftarrow m_{s_i,i} + \rho_{s_i}
\end{aligned}
\end{equation}
    We then set $i \leftarrow i+1$ and continue until $i = \binom{d}{k}.$ 
    \item{Step 2:} Set $count \leftarrow count+1$ and return to step 1. We also set the convergence condition as the variance of the tensor elements projected by each of its dimensions is smaller than a very small positive $\epsilon$.
    
\end{enumerate}
The time complexity of this algorithm is $ O\left(\sum\limits_{i=1}^{\binom{d}{k}} \sum\limits_{s_i = 1}^{|(i)|}\phi(s_i,i) \right)\equiv O\left(\binom{d}{k}\norm{\upsigma(A)} \right).$ This complexity is unaffected by the final step of converting our solution over logarithms to our desired solution as $A(\vec{\alpha}) = \exp(a(\Vec{\alpha}))$, which is proportional to the number of nonzero elements of $A$.

\section{Discussion}

In this paper we have generalized the RZ algorithm \cite{RZalgorithm} as specialized for the unique canonical scaling of matrices \cite{UC18,Aktar19,UC15} for the canonical scaling of arbitrary $k$-dimensional subtensors of a given $n$-dimensional tensor. This algorithm provides a means for generalizing scale-invariant methods presently applied to matrices to higher-dimensional structures, e.g., extending the scaling-based image interpolation method of \cite{Aktar19} to the interpolation of a video sequence treated as a 3d array/tensor.

\bibliographystyle{plain}

\clearpage

\appendix
\section{Appendix 1: Proof of the KT construction}
Here we establish formulas to establish the equivalency of \textbf{Program I} and \textbf{Program II} in section 3.1.3. We define the unfolding tensor column vector based on formula \ref{gettheorderforC}. For equation $Cx=b$, all possible elements of the set $s = \{s_{s_i,i}\}$ maps onto $b$ by giving an ordering for $(s_i, i)$ when putting all such pairs in a column vector. A function $P$ with an ordering for $(s_i, i)$ is defined as
\begin{equation}
    P(s_i,i) = \sum_{i' = 0}^{i-1} |(i')| + s_i,
\end{equation}
with the assumption $|(0)| = 0$. Then all elements $s_{s_i,i}$ maps onto $b$ and all elements $m_{s_i,i}$ maps onto $\omega$ as
\begin{equation}
    b\left(P(s_i,i)\right) = s_{s_i,i}, \quad 
    \omega\left(P(s_i,i)\right) = m_{s_i,i}. 
\end{equation}
Note that the tensor $x$ agrees with tensor $a$ in terms of nonzero elements. With a slight abuse of notation, we transform tensor $x$ to a vector $x$ by unfolding vector notations from definition \ref{gettheorderforC}. 
Matrix $C$ becomes the medium for both \textbf{Program I} and \textbf{Program II}, as $\omega^T C$ at position $J(\Vec{\alpha})$ picks the corresponding scaling elements and $Cx$ at position $P(s_i,i)$ picks nonzero elements of a subtensor $A_{s_i,i}$. The matrix $C$ is defined with its column size equals $N = \prod_{i=1}^d n_i$ and its row size equals $\sum_{i=1}^{\binom{d}{k}}|(i)|.$

\begin{equation}          
 C_{P(s_i,i) J\left(\Vec{\alpha}\right)}  \equiv  \left\{ \begin{aligned}
   &1 \quad \text{if} \quad   \Vec{\alpha} \in \upsigma_A(s_i,i) \\  
&0 \quad \text{  otherwise.}
\end{aligned}\right.
\end{equation}
Then
\begin{equation}
(Cx)_{P(s_i,i)} = \sum\limits_{J(\Vec{\alpha})} C_{P(s_i,i), J(\Vec{\alpha})} x\left({J(\Vec{\alpha})}\right) = \sum\limits_{\Vec{\alpha} \in \upsigma_{A}(s_i,i)} x(\Vec{\alpha}).     
\end{equation}
Thus, $Cx = b$ is equivalent to
\begin{equation}
  \sum\limits_{\Vec{\alpha} \in \upsigma_{A}(s_i,i)} x(\Vec{\alpha}) = s_{s_i, i}.  
\end{equation}
Likewise, equation $\omega^T C$ satisfies 
\begin{equation}
    [\omega^T C]_{J(\Vec{\alpha})} = \sum\limits_{P(s_i,i)}  \omega\left( P(s_i,i)\right) C_{P(s_i,i), J(\Vec{\alpha})} = \sum\limits_{(s_i, i) \in \upsigma_A(\Vec{\alpha})} m_{s_i, i}  
\end{equation}

\section{Appendix 2: Characterization Proof}
We establish the proof for \textbf{Theorem  \ref{characterize}}. Using the proof from Section 3 of \cite{RZalgorithm} is sufficient to establish from conditions $1$ to condition $5$. Since condition $6$ is equivalent to condition $7$, we only show the equivalency of condition $7$ with former conditions. Let's consider 
$z$  such that $z\left(P(s_i,i)\right) = \mu_{s_i,i} $ and $Cx = b$ could assert that:
\begin{equation}
    z^{T} C \left[J(\Vec{\alpha})\right] = \sum\limits_{(s_i,i) \in \sigma_A(\Vec{\alpha})} \mu_{s_i,i}  = 0
\end{equation}
Then by the following lemma,

\begin{lemma}
If $Cx = b$ and there exists a vector $\lambda$ with $\lambda^T C = 0$ then $\lambda ^T b = 0.$
\end{lemma}{}
We obtain condition $7$ and the theorem is fully proved, followed by
\begin{equation}
\Vec{z}^T\Vec{b} =  
    \sum\limits_{i=1}^{\binom{d}{k}} \sum\limits_{s_i = 1}^{|(i)|} s_{s_i,i}\cdot\mu_{s_i,i} = 0. 
\end{equation}

\end{document}